\documentclass[a4paper,12pt]{amsart}
\usepackage[utf8]{inputenc}
\usepackage[english]{babel} 

\usepackage[T1]{fontenc}
\input{diagrams.tex}
\diagramstyle[scriptlabels,height=8mm,width=8mm]
\usepackage{amsthm}
\usepackage{amsmath}
\usepackage{amsfonts}
\usepackage{mathtools}
\usepackage[colorlinks]{hyperref}
\usepackage{enumerate}
\usepackage{amssymb}
\usepackage{tikz,pgf}
\usepackage{graphicx}
\usetikzlibrary{matrix,arrows,decorations.pathmorphing}
\usetikzlibrary{calc}
\makeatletter
\newbox\dottedarrow@box
\setbox\dottedarrow@box\hbox
  {%
    \begin{tikzpicture}
      \draw[dotted,->] (0,0) -- (1.5em,0);
    \end{tikzpicture}%
  }
\newcommand*\dottedarrow
  {\relax\ifmmode\expandafter\dottedarrow@m\else\expandafter\dottedarrow@t\fi}
\newcommand*\dottedarrow@t[1][1.5em]
  {\resizebox{#1}{!}{\raisebox{.5ex}{\usebox\dottedarrow@box}}}
\newcommand*\dottedarrow@m[1][]
  {%
    \if\relax\detokenize{#1}\relax
      \mathchoice
        {\dottedarrow@t}
        {\dottedarrow@t}
        {\dottedarrow@t[1.1em]}
        {\dottedarrow@t[0.9em]}%
    \else
      \dottedarrow@t[#1]%
    \fi
  }
\makeatother

\title[Algebraic Gromov ellipticity]{
Algebraic Gromov ellipticity: a brief survey}
{}
\author{Mikhail Zaidenberg}

\address{Institut Fourier, UMR 5582, 
Laboratoire de Math\'ematiques,\newline\indent
Universit\'e Grenoble Alpes, 
CS 40700, 38058 Grenoble cedex 9, France
}
\email{mikhail.zaidenberg@univ-grenoble-alpes.fr}
\date{}

\newtheorem{thm}{Theorem}[section]
\newtheorem{cor}[thm]{Corollary}
\newtheorem{lem}[thm]{Lemma}
\newtheorem{prop}[thm]{Proposition}
\theoremstyle{definition}
\newtheorem{defi}[thm]{Definition}
\newtheorem{defis}[thm]{Definitions}
\newtheorem{conj}[thm]{Conjecture}
\newtheorem{conv}[thm]{Convention}
\newtheorem{nota}[thm]{Notation}
\newtheorem{rem}[thm]{Remark}
\newtheorem{rems}[thm]{Remarks}
\newtheorem{exa}[thm]{Example}
\newtheorem{exas}[thm]{Examples}
\newtheorem{prob}[thm]{Problem}
\newtheorem{probs}[thm]{Problems}
\newtheorem{ques}[thm]{Question}
\newtheorem{sett}[thm]{Setting}
\newtheorem{sit}[thm]{}

\newcommand{\Aut}{ \operatorname{{\rm Aut}}}
\newcommand{\SAut}{ \operatorname{{\rm SAut}}}
\newcommand{\End}{ \operatorname{{\rm End}}}

\newcommand{\PGL}{\operatorname{{\rm PGL}}}

\renewcommand{\epsilon}{\varepsilon}

\def\and{\quad\mbox{and}\quad}

\newcommand{\C}{\ensuremath{\mathbb{C}}}

\newcommand{\GG}{\ensuremath{\mathbb{G}}}

\newcommand{\A}{\ensuremath{\mathbb{A}}}

\newcommand{\kk}[1]{\bk^{[#1]}}

\newcommand{\tE}{{\tilde E}}

\newcommand{\tp}{{\tilde p}}

\newcommand{\ts}{{\tilde s}}

\def\fg{{\mathfrak g}}

\renewcommand{\rho}{\varrho}

\def\bals#1\eals{\begin{align*}#1\end{align*}}
\def\bal#1\eal{\begin{align}#1\end{align}}

\def\SAut{\mathop{\rm SAut}}

\def\kk{{\mathbb K}}
\def\A{{\mathbb A}}

\def\TT{{\mathbb T}}
\def\CC{{\mathbb C}}

\def\PP{{\mathbb P}}

\renewcommand{\phi}{\varphi}

\newcommand{\bnum}{\begin{enumerate}}
\newcommand{\enum}{\end{enumerate}}

\newcommand{\brem}{\begin{rem}}
\newcommand{\brems}{\begin{rems}}
\newcommand{\erem}{\end{rem}}
\newcommand{\erems}{\end{rems}}
\newcommand{\bprob}{\begin{prob}}
\newcommand{\eprob}{\end{prob}}
\newcommand{\bprobs}{\begin{probs}}
\newcommand{\eprobs}{\end{probs}}
\newcommand{\bques}{\begin{ques}}
\newcommand{\eques}{\end{ques}}
\newcommand{\bexa}{\begin{exa}}
\newcommand{\bexas}{\begin{exas}}
\newcommand{\eexa}{\end{exa}}
\newcommand{\eexas}{\end{exas}}
\newcommand{\bdefi}{\begin{defi}}
\newcommand{\edefi}{\end{defi}}
\newcommand{\bdefis}{\begin{defis}}
\newcommand{\edefis}{\end{defis}}
\newcommand{\bcor}{\begin{cor}}
\newcommand{\ecor}{\end{cor}}
\newcommand{\blem}{\begin{lem}}
\newcommand{\elem}{\end{lem}}
\newcommand{\bconv}{\begin{conv}}
\newcommand{\econv}{\end{conv}}
\newcommand{\bconj}{\begin{conj}}
\newcommand{\econj}{\end{conj}}
\newcommand{\bprop}{\begin{prop}}
\newcommand{\eprop}{\end{prop}}
\newcommand{\bthm}{\begin{thm}}
\newcommand{\ethm}{\end{thm}}
\newcommand{\bnota}{\begin{nota}}
\newcommand{\enota}{\end{nota}}
\newcommand{\bsit}{\begin{sit}}
\newcommand{\esit}{\end{sit}}
\newcommand{\be}{\begin{equation}}
\newcommand{\ee}{\end{equation}}
\newcommand{\bproof}{\begin{proof}}
\newcommand{\eproof}{\end{proof}}
\newcommand{\bsett}{\begin{sett}}
\newcommand{\esett}{\end{sett}}
\def\ba{\begin{array}}
\def\ea{\end{array}}

\makeatletter
\def\blfootnote{\xdef\@thefnmark{}\@footnotetext}
\makeatother
\begin{document}
\begin{abstract} We survey on 
algebraically elliptic varieties in the sense of Gromov.
\end{abstract}
\maketitle
 \mbox{\hspace{20pt}}

\date{}

{\footnotesize \tableofcontents}
\bigskip

\section{Introduction}\label{sec:intro}
Gromov ellipticity is often considered to be
a chapter of complex analysis. 
However, in his foundational paper \cite{Gro89}
Gromov also presented and studied 
an algebraic counterpart of this notion.
In this survey, we mainly focus on the algebraic part.
All algebraic varieties considered below
are defined over an algebraically closed field 
$\kk$ with
$\rm{char}\, \kk=0$. 
In what follows by \emph{Gromov spray} and 
\emph{Gromov ellipticity} 
we mean 
\emph{algebraic Gromov spray} and 
\emph{algebraic Gromov ellipticity}; 
otherwise we talk about 
\emph{analytic Gromov spray} and 
\emph{analytic Gromov ellipticity}.
We address the comprehensive monograph \cite{For17}
and the survey article \cite{For23} by F.~Forstneri\v{c}
for a thorough introduction to Gromov ellipticity, 
especially for its complex analytic counterpart. 
\subsection{Prehistory: the Oka-Grauert principle}
The origin of Gromov ellipticity in complex analysis lies in 
the following Oka-Grauert principle; see 
\cite[Sec.~5.1]{For17} for a historical account. 
\begin{thm}[\rm{see Oka \cite{Oka39}, Frenkel \cite{Fre57}, 
Grauert \cite{Gra57}-\cite{Gra58}, A.~Cartan \cite{Car58}, 
Ramspott \cite{Ram65}, Henkin--Leiterer \cite{HeLa98}}] $\,$

\noindent Given a complex Lie group $G$, 
the classifications of principal 
$G$-fiber bundles over a Stein complex manifold $S$ 
in topological and holomorphic categories coincide. 
The same holds for the classes up to homotopy of sections 
of the associated fiber bundles with 
$G$-homogeneous fibers.
\end{thm}
Among important proprerties of a complex Lie group $G$ 
linked to the Oka-Grauert principle,
we distinguish the following. 
\begin{thm}[\rm{Grauert \cite{Gra57}, 
Gromov \cite[Sec.~1.4D$'$]{Gro89}}] $\,$

\noindent For a complex Lie group $G$ 
and a complex Stein manifold $S$, 
the following hold.
\begin{enumerate}
\item[{\bf  (A)}]  
Every continuous map $S\to G$ 
is homotopic to 
a holomorphic map $S\to G$.
\item[{\bf (B)}] 
Every holomorphic map $\bar D\to G$, 
where $D\subset \CC^n$ 
is a bounded convex domain, 
can be approximated by  holomorphic 
maps $\CC^n\to G$ 
uniformly on $\bar D$. 
\end{enumerate}
\end{thm}
\subsection{Analytic Gromov ellipticity} 
Loosely speaking, 
in his paper \cite{Gro89} Gromov
answers the following  question.
\begin{ques} \emph{What do you need to know 
about a complex manifold $X$ 
to be sure that the analogues of {\bf (A)} 
and {\bf (B)} hold for holomorphic maps $S\to X$? }
\end{ques}

The following answer is a manifestation of Gromov's 
$h$-principle for complex manifolds
(see subsection \ref{ss:h-principle} below):

\medskip

\noindent \emph{Every analytically
elliptic complex manifold $X$~verifies  
analogues of~{\bf(A)} and~{\bf (B)}. }
\medskip

In the following two subsections, we give 
the definitions of Gromov spray and 
Gromov ellipticity. 
\subsection{Gromov spray}

\begin{defis} [Gromov \cite{Gro89}] 
Let $X$  be  a complex manifold. 
A \emph{spray} on $X$ is a triple $(E,p,s)$, where
\begin{enumerate}
\item[$\bullet$] $p\colon E\to X$ is 
a holomorphic vector bundle on $X$
with zero section $Z$, and 
\item[$\bullet$] $s\colon E\to X$ is 
a holomorphic map such that 
$s|_Z=p|_Z$, i.e.
$s(0_x)=x$ for all $x\in X$, where $0_x$ is the origin  
of the vector space $E_x=p^{-1}(x)$. 
\end{enumerate}

A spray $(E,p,s)$ on $X$ is called \emph{dominating} 
if for every $x\in X$ the differential  
$ds\colon T_{0_x}E_x\to T_xX$ is onto.
\end{defis}
\subsection{Ellipticity and subellipticity}
\begin{defi} [\rm{cf. Gromov 
\cite[Sec.~0.5 and 3.5B]{Gro89}}]$\,$
\begin{enumerate}
\item $X$ is called \emph{analytically elliptic} 
if it admits a dominating spray $(E,p,s)$.
\item $X$ is called \emph{locally analytically elliptic} 
if it admits an open covering $(U_i)_i$ 
with dominating sprays $(E_i,p_i,s_i)$ on $U_i$ 
 with values in $X$, i.e.
$s_i$  is a holomorphic map $U_i\to X$. 
\end{enumerate}
\end{defi}
\begin{defi} [\rm{Forstneri\v{c} 
\cite[Definition 2]{For02}}]
$X$ is called \emph{analytically subelliptic} 
if it admits a dominating family 
of sprays $(E_i,p_i,s_i)$ on $X$, i.e.
for each $x\in X$ we have
\begin{equation}\label{eq:subelliptic}
T_xX={\rm span} \left(\bigcup_i 
ds_i(T_{0_x}E_{i,x})\right).
\end{equation}
\end{defi} 
\subsection{Gromov's h-principle for 
complex manifolds}\label{ss:h-principle}
\begin{defi}[\rm{Gromov \cite[Sec.~0]{Gro89}}] 
Let $X$ and $Y$ be complex manifolds. 
We say that holomorphic maps $Y \to X$ 
satisfy the $h$-principle ($h$ for homotopy) 
if every continuous map $Y \to X$ is homotopic 
to a holomorphic map.
\end{defi}
\begin{thm} [\rm{Gromov \cite[Sec.~1.4D$'$]{Gro89}, 
Forstneri\v{c} \cite[Theorem 1.2]{For06b}}]
\label{thm:A-B-prim}
\label{thm:CAP} $\,$

\noindent Let $S$ be a Stein complex space and $X$ be
a complex manifold. If $X$ is analytically elliptic, 
then the following hold.
\begin{enumerate} \item[{\bf (A$'$)}] 
 Every continuous map $S\to X$
  is homotopic to a holomorphic one;
\item[] the same holds for sections of 
holomorphic fiber bundles over $S$ with fiber $X$.
 \item[{\bf (B$'$)}] Every holomorphic map 
 $f \colon \bar D \to X$, where $D \subset \CC^n$ 
is a bounded convex domain,  can be 
uniformly on $\bar D$ approximated 
by holomorphic maps $\CC^n \to X$. 
\end{enumerate} 
\end{thm}
A complex manifold $X$ verifying  ({\bf B$'$}) 
is called an \emph{Oka manifold}, 
see \cite{For06a}, \cite{For09} 
and \cite[Definition~5.4.1]{For17}. 
An Oka manifold verifies a stronger 
\emph{convex approximation property} 
(CAP, for short), see \cite[Theorem~5.4.4]{For17}.
Any analytically elliptic complex manifold 
is Oka, it verifies the CAP and a condition  Ell$_1$ 
of Gromov, see  \cite[Proposition~8.8.11]{For17}. 
Actually,
the Oka property is equivalent 
to the condition Ell$_1$
(Kusakabe \cite[Theorem~1.3]{Kus21b})
and does not imply analytic ellipticity, in general,
see Kusakabe \cite{Kus20b} 
and \cite[Corollary~1.5.]{Kus24}
for corresponding examples. 
\subsection{Algebraic Gromov ellipticity}
As we have already mentioned,
Gromov \cite[Sec.~3.5A]{Gro89} also introduced 
the notions of 
algebraic spray, algebraic ellipticity, etc.,
where complex manifolds are replaced 
by smooth algebraic varieties (i.e. algebraic manifolds) 
defined over $\kk$, 
holomorphic vector bundles by algebraic vector bundles and 
holomorphic maps by regular maps. 
In the algebraic category, we have
 the following equivalences.
\begin{thm} [\rm{Gromov \cite[Sec.~3.5.B$'$]{Gro89}, 
Kaliman--Zaidenberg \cite[Theorem~1.1]{KaZa24a}}]
\label{thm:equiv} $\,$

\noindent For a smooth algebraic variety $X$, 
the following are equivalent:
\begin{enumerate}
\item $X$ is algebraically Gromov elliptic;
\item $X$ is locally algebraically Gromov elliptic;
\item $X$ is algebraically subelliptic.
\end{enumerate}
\end{thm}
An analogue of (1)$\Leftrightarrow$(3) 
in the analytic category is known to hold
for complex Stein manifolds $X$, see \cite[Lemma~2.2]{For02}. 
 
 The proof of (2)$\Rightarrow$(1) uses the following 
 Localization Lemma of Gromov \cite[Sec.~3.5.B]{Gro89}; 
 see also \cite[Propositions~6.4.1--6.4.2]{For17} 
 and \cite[Proposition 8.1]{KaZa23b}.
 \begin{lem}\label{lem:loc}
 Let $X$ be a smooth algebraic variety, 
 $D$ be a reduced 
 effective divisor on $X$ and $(E,p,s)$ 
 be a Gromov spray on 
 $U = X \setminus {\rm supp}(D)$ with values in $X$ 
 such that $p\colon E \to U$ is a trivial vector bundle. 
 Then there exists a spray $(\tE, \tp, \ts)$ 
 on $X$ whose restriction to $U$ is isomorphic
to $(E,p,s)$. In particular, if $(E,p,s)$ 
is dominating on $U$, then so is 
$(\tE, \tp, \ts)|_U$. 
 \end{lem}
Yet another important ingredient  in
the proofs of (2)$\Rightarrow$(1) 
and (3)$\Rightarrow$(1) 
is the composition of Gromov sprays, 
see Gromov \cite[Sec.~1.3.B]{Gro89}. 
Given two sprays $(E_1,p_2,s_1)$ and 
$(E_2,p_2,s_2)$
on $X$, we consider their composition $(E,p,s)$ where 
\begin{itemize}
\item[$\bullet$] $E= \{(e_1, e_2) \in E_1\times E_2
\vert  s_1(e_1) = p_2(e_2)\}$, 
and so
${\rm pr}_1|_E\colon E=s_1^*E_2\to E_1$ 
is the induced vector bundle; 
\item[$\bullet$]
$p=p_1\circ {\rm pr}_1\colon E\to X$;
\item[$\bullet$] 
$s=s_2\circ {\rm pr}_2\colon E\to X$. 
\end{itemize}
In general, $p\colon E\to X$
is not a projection of a vector bundle. 
However, $p$ happens to be 
such a projection provided that 
$p_2\colon E_2\to X$ is a line bundle, 
see \cite[Proposition~2.1]{KaZa24a}. 

If $X$ is algebraically subelliptic, 
then we can find 
$m\ge\dim(X)$ rank 1 sprays
$(L_i,p_i,s_i)$ on $X$ which satisfy 
an analogue of \eqref{eq:subelliptic}.
The iterated composition 
of these sprays yields 
a dominating Gromov spray $(E,p,s)$ 
on $X$ of rank $m$ and 
provides the ellipticity of $X$.

It is easily seen that the product of 
Gromov elliptic  smooth
algebraic varieties is Gromov elliptic. 
The converse is also true.
\begin{lem}[\rm{L\'arusson's Lemma, 
see e.g. \cite[Lemma~3.6]{AKZ24}}]\label{lem:LL} $\,$
 If the product $X_1\times X_2$ 
of two smooth algebraic varieties is Gromov elliptic, 
then $X_1$ and $X_2$ are Gromov elliptic.
\end{lem}
\begin{lem}[\rm{Forstneri\v{c}
 \cite[Proposition~6.4.10]{For17},
 L\'arusson--Truong 
\cite[Theorem~1 and Remark 2(a)]{LT19}}] $\,$
\label{lem:cover}

\noindent Let $\tilde X\to X$
be a proper \'etale morphism of smooth 
complex algebraic varieties. If $X$ 
is Gromov elliptic, then also $\tilde X$ is. 
\end{lem}
\section{Flexible varieties}
\subsection{Flexibility versus Gromov ellipticity}
\begin{exa}\label{ex:Ga} 
Let $\mathbb{G}_a=(\kk,+)$ 
be the additive group of the field $\kk$, let $X$ be 
a smooth algebraic $\mathbb{G}_a$-variety
and $s\colon {\mathbb G}_a\times X\to X$ 
be the action morphism. 
 Consider the trivial line bundle  
 $\,\,p\colon L={\mathbb G}_a\times X\to X$ of rank 1, 
 where $p$ is the second projection. 
 Then $(L,p,s)$ is a rank 1 spray on $X$ dominating 
 in directions of the one dimensional
 ${\mathbb G}_a$-orbits $s(E_x)={\mathbb G}_a.x$. 
 The latter means
 that ${\rm rank}\, (d s|_{T_{0_x}E_x})=1$ 
 provided~$\dim({\mathbb G}_a.x)=1$. 
 \end{exa}
\begin{defi} A smooth 
quasiaffine algebraic variety $X$  is called  \emph{flexible}
if there exists a collection $U_1,\ldots,U_k$ of 
${\mathbb G}_a$-subgroups of $\Aut(X)$ such that
for each $x\in X$ the velocity vectors 
of $U_1,\ldots,U_k$ at $x$ 
span $T_xX$.
It is said to be \emph{locally flexible} is 
$X$ admits a Zariski open covering by 
flexible quasiaffine varieties.
\end{defi}
Given a flexible algebraic manifold $X$, 
consider rank 1 sprays 
$(L_i,p_i,s_i)$ associated with the 
${\mathbb G}_a$-subgroups $U_1,\ldots,U_k$ as above, 
see Example \ref{ex:Ga}. The
composition of these sprays
provides the subellipticity of $X$. 
Together with Theorem \ref{thm:equiv} 
this leads to the following result, see 
Gromov \cite[Sec.~0.5.B]{Gro89} and
also \cite[Appendix]{AFKKZ13},
\cite[Proposition~5.6.22(C)]{For17} 
and \cite[Theorem~3.1]{KKT18}.
 \begin{prop}\label{thm:loc-flexible}
A locally flexible smooth  algebraic 
variety is Gromov elliptic. 
 \end{prop}
 \begin{defi}\label{def:loc-stab-flex}
One says that a smooth algebraic variety $X$ is 
\emph{stably flexible}
if $X\times\A^k$ is flexible for some $k\ge 0$, and 
\emph{locally stably flexible} if $X$ 
admits an open covering
by stably flexible affine charts. 
\end{defi}
If $X$ admits a covering
by copies of $\A^n$, then it is certainly 
locally flexible. Using Theorem \ref{thm:equiv}, 
Proposition \ref{thm:loc-flexible}
and L\'arusson's Lemma \ref{lem:LL} one 
deduces the following strengthening of 
Proposition \ref{thm:loc-flexible}.
\begin{thm}[\rm{cf. \cite[Corollary~3.2]{KKT18}}]
\label{thm:loc-stab-flex} 
Every locally stably flexible
smooth algebraic variety is Gromov elliptic. 
\end{thm}
Let us give examples of non-flexible, 
but stably flexible affine surfaces.
\begin{exa}[\rm{\cite[Example~0.4]{KKT18}}] 
\label{ex:Dan}
Consider the smooth 
affine Danielewski surfaces $S_k$ 
given in $\CC^3$ by equations 
$\{x^ky-z^2+1=0\}$. 
The surface $F_1$ is flexible,
while $F_k$ with  $k>1$ 
is not. This follows e.g. from the description of 
the automorphism groups $\Aut(F_k)$, 
see Makar-Limanov \cite{ML01}. 
For every $k\ge 2$  we have $F_k\times \C\cong  F_1\times \C$ 
(Danielewski \cite[Theorem 1]{Dan89}). 
Since $F_1\times \C$ is flexible, 
$F_k$ with  $k>1$ is stably flexible, 
while being non-flexible. 
Due to Theorem \ref{thm:loc-stab-flex}, 
$F_k$ is Gromov elliptic 
for every $k\ge 1$.
\end{exa}
Another characteristic property of flexible varieties 
are their homogeneity properties,
see \cite[Theorem~1.1]{AFKKZ13},
\cite[Theorem~1.1]{FKZ16}, and  
\cite[Theorem~11]{Arz18}.
\begin{thm}\label{thm:high-transitivity}
 Let $X$ be a smooth quasiaffine variety of dimension $\ge 2$ 
 and $\SAut(X)$ be the subgroup of $\Aut(X)$ generated by 
 all $\GG_{\mathrm a}$-subgroups of $\Aut(X)$. 
 Then $X$ 
is flexible if and only if $\SAut(X)$ acts transitively on $X$,
if and only if it acts $m$-transitively 
on $X$ for any $m\ge 1$.
\end{thm}
\begin{defi}[\rm{Bogomolov-Karzhemanov-Kuyumzhiyan \cite[Definition~1.2]{BKK13}}] $\,$

\noindent One says that an algebraic variety $X$ is \emph{birationally stably flexible}
if the field  extension $\kk(X)(y_1,\ldots,y_n)$ admits a flexible model.  
\end{defi}
See \cite[Theorems~2.1 and~2.2]{BKK13} for criteria of the birational stable flexibility.
Clearly, any stably rational variety is birationally stably flexible. On the other hand,
a birationally stably flexible variety is unirational. There is the following conjecture. 
\begin{conj}[\rm{\cite[Conjecture~1.4]{BKK13}}]
\emph{Any unirational algebraic variety $X$ is birationally stably flexible.}
\end{conj}
\subsection{Examples of flexible varieties and of
Gromov elliptic varieties}\label{ss:flex-ex}
\begin{enumerate}
\item
 Let $G$ be a connected complex Lie group 
 with Lie algebra $\fg$, let $Y=G/H$ be 
 a homogeneous manifold of $G$ and 
 $\exp \colon \fg \to G$ be the exponential
map. Then the map 
$Y \times \fg\to Y$, $(y,v)\mapsto \exp(v)y$ is 
a dominating analytic spray on $Y$, 
see \cite[Sec.~3.1]{For23}. Hence,
the homogeneous space $Y$ is analytically elliptic. 

In particular,
$\CC^*=\CC\setminus\{0\}$ 
is elliptic in analytic sense, 
but is not elliptic in algebraic sense. 
The same holds for any smooth projective curve of genus 1.

\item If $G$ is a semisimple linear algebraic group,
then $G/H$ is flexible, see \cite[Proposition~5.4]{AFKKZ13}, 
and so Gromov elliptic. 
In particular, every flag variety $G/P$, where $P\subset G$ 
is a parabolic subgroup, is Gromov elliptic.

\item The affine space $\A^n$ is Gromov elliptic. 

\item Let $X$ be an algebraic variety. If $X$  admits 
a Zariski open covering $(U_i)$ 
where $U_i\simeq\A^n$,
then $X$ is locally Gromov elliptic, 
hence Gromov elliptic by Theorem \ref{thm:equiv}. 

\item  Every  smooth complete spherical 
variety $X$ 
admits an open covering by affine spaces 
(Brion--Luna--Vust \cite[Sec.~1.5, Corollaire]{BLV86}).
By (4), $X$ is Gromov elliptic.
In particular, every smooth complete toric 
variety  is Gromov elliptic. 
Moreover, a smooth toric variety 
with no torus factor
is covered by affine spaces, 
hence is Gromov elliptic, cf. \cite[Remark~4.7]{AKZ24}.
 
\item
A smooth hypersurface $X$ in $\A^{n+2}$ 
given by equation
\[uv-p(x_1,\ldots,x_n)=0,\]  
where $p\in\kk[x_1,\ldots,x_n]$ is 
a nonconstant polynomial,
is flexible, see \cite[Theorem~5.1]{KaZa99}  
and \cite[Theorem~0.1]{AFKKZ13}, or alternatively 
\cite[Theorem~0.2(3)]{AKZ12}. 
So, $X$ is Gromov elliptic (cf. Example \ref{ex:Dan}). 

\item Every smooth complete rational 
surface $S$ admits 
a covering by copies of $\A^2$. 
By (4), $S$ is Gromov elliptic. 
 
\item  A flexible smooth quasiaffine algebraic 
variety remains flexible after
deleting a closed subvariety of codimension $\ge 2$ 
(Flenner-Kaliman-Zaidenberg \cite[Theorem~0.1]{FKZ16}). 
Likewise, a locally flexible algebraic variety
remains locally flexible after
deleting a closed subvariety of codimension $\ge 2$.

\item Let $X$ be a smooth algebraic variety 
covered by open charts isomorphic to $\A^n$
and $A\subset X$ be a closed algebraic 
subvariety  of codimension $\ge 2$. 
Then $X\setminus A$ is  Gromov elliptic 
(see Forstneri\v{c}
\cite[Proposition~6.4.5]{For17}). 
Actually, this follows immediately from (8) 
due to Theorems \ref{thm:equiv} and 
\ref{thm:loc-flexible}. 
Cf. also Gromov \cite[Sec. 3.5C]{Gro89} and 
Kusakabe \cite{Kus21b} for stronger results.

\item The configuration spaces of a flexible 
quasiaffine variety $X$
of dimension $\dim(X)\ge 2$ are flexible 
(Kusakabe \cite[Proposition~3.4]{Kus24}).
\end{enumerate}
 For more examples of flexible 
 varieties, see \cite[Sec.~4.1]{CPPZ21}. 
See also Theorem \ref{thm:flex-cones} below
for examples of non-flexible Gromov 
elliptic quasiaffine varieties.

For the following classes of birationally 
stably flexible varieties, 
see \cite[Sections 3.1-3.4]{BKK13}. 
\begin{exas} 
1.  Let $G \subset \PGL(n+1,\kk)$ be a finite subgroup.
Then the quotient 
$\PP^n/G$ is birationally stably flexible.

2.
Let $X$ be an algebraic variety. 
Assume 
that $X$ carries a collection of distinct 
birational structures of $\PP^{m_i}$-bundles, 
$\pi_i \colon X \to S_i$ such that the tangent 
spaces of generic fibers of $\pi_i$ 
span the tangent space of $X$ at the generic point. 
Then $X$ is birationally stably flexible.

3. Every smooth  cubic hypersurface $X\subset \PP^{n+1}$, $n \ge 2$, 
 is birationally stably flexible. The same holds for quartic hypersurfaces 
$ X \subset \PP^{n+1}$, $n \ge 3$ that have a line of double singularities,
and for smooth complete intersections of three quadrics in $\PP^6$.
\end{exas}
\section{Properties of Gromov elliptic varieties}
\subsection{Approximation results}
In the algebraic category we have 
the following analogue of the Oka property {\bf (B$'$)}. 
\begin{thm}[\rm{Forstneri\v{c} \cite[Corollary~6.15.2]
{For17}}]\label{thm:approx} $\,$

\noindent Let $X$ be a smooth complex 
algebraic variety and $D$ 
be a bounded convex domain in $\CC^n$. 
If $X$ is Gromov elliptic, then any 
holomorphic map $\bar D\to X$ 
can be approximated by morphisms 
$\CC^n \to X$ uniformly on $\bar D$. 
\end{thm} 
\begin{rem}\label{rem:aCAP}
The approximation  property in the conclusion 
of Theorem \ref{thm:approx}
is called 
the \emph{algebraic convex approximation property}, 
abbreviated as \emph{aCAP}. 
Thus, according to Theorem \ref{thm:approx}
Gromov ellipticity implies aCAP. 
We don't know if the converse is true. 
\end{rem}

According to Forstneri\v{c} \cite[Theorem~1.1 
and Corollary~1.2]{For06a} and 
L\'arusson-Truong \cite[Theorem 1]{LT19}, 
the algebraic Gromov ellipticity is 
equivalent to two other important 
properties. 
\begin{defi}[\rm{\cite[Definition~1.1]{Kus21a}}]
\label{def:aOka}
A smooth algebraic variety $X$ over $\kk$ is called
 \emph{algebraically Oka} 
(abbreviated as  \emph{aOka}, or $aEll_1$ after Gromov)
if for each morphism $f \colon Y \to X$, where $Y$  is
an affine algebraic manifold, 
there exists a morphism 
$F\colon Y\times \A^N \to X$ such that
$F|_{Y\times\{0\}}= f$ and 
$F|_{\{y\}\times \A^N}\colon \A^N \to X$ 
is a submersion at $0\in\A^N$ for each $y\in Y$, 
i.e. \[dF|_{T_{(y,0)}(\{y\}\times\CC^m)}\colon 
T_{(y,0)}(\{y\}\times\A_{\CC}^m)\to T_{f(y)}X\] 
is onto for every $y\in Y$.
\end{defi}
\begin{thm}[\rm{L\'arusson-Truong 
\cite[Theorem 1]{LT19}}]\label{thm:LT}  $\,$

\noindent For a smooth complex algebraic 
variety $X$, 
the following are equivalent.
\begin{enumerate}
\item $X$ is Gromov elliptic;

\item $X$ is aOka;

\item given a morphism
$f_0 \colon Y\to X$ from 
a complex affine variety $Y$, 
a holomorphically convex compact subset 
$K\subset Y$, 
a subset $U\supset K$ of $X$ open in the 
complex topology of $X$, 
and a homotopy of holomorphic maps 
$f_t\colon U \to X$, 
$t \in [0, 1]$, there is a morphism
$F \colon Y\times \A^1_{\CC}\to X$ 
with $F(\cdot,0)=f_0$ 
and $F(\cdot,t)$ as close to $f_t$ as desired, 
uniformly on $K$.
\end{enumerate}
\end{thm}
In the analytic setup, the properties 
analogous  to (1)--(3) of 
Theorem \ref{thm:LT} are 
known to be equivalent provided that
 the complex manifold $X$ is Stein, 
 see e.g. \cite[Sec.~5.5]{For17}, 
\cite{LT19} and the references therein. 
\begin{rem}\label{rem:Lefschetz} 
In fact, the equivalence (1)$\Leftrightarrow$(2) 
of Theorem \ref{thm:LT}
holds for algebraic varieties defined 
over a general algebraically closed field $\kk$
of characteristic zero. This can be shown by 
inspecting the arguments in the proof. 
For the reader's convenience, we provide 
a  proof
of the equivalence (1)$\Leftrightarrow$(2) 
in this generality.
\end{rem}
\begin{proof} 
 (1)$\Rightarrow$(2). Assume that $X$ is Gromov elliptic, 
 and let $(E,p,s)$ be a dominating spray over $X$. 
 Let $Y$ be an affine variety and 
 $f\colon Y\to X$ be a morphism. 
 Consider the pullback $(\tE,\tp,\ts)=f^*(E,p,s)$, 
 where $\tp\colon\tE\to Y$ is the induced vector bundle
 and $\ts=s\circ f\colon \tE\to X$. Since $Y$ is affine, 
 the vector bundle $\tp\colon\tE\to Y$
 is globally generated, 
 see Serre \cite[Sec.~45, P.~238, Corollaire~1]{Ser55}. 
 Let $\eta_1,\ldots,\eta_m$ 
 be global sections of $\tp\colon\tE\to Y$ 
 that span every fiber $\tp^{-1}(y)$, $y\in Y$.
 Then the morphism of vector bundles 
 \[\phi\colon Y\times\A^m\to_Y \tE,\quad 
 (y, (a_1,\ldots,a_m))\mapsto \sum_{i=1}^m a_i\eta_i(y)\]
 identical on the base $Y$ yields a fiberwise linear surjection.
 It is easily seen that the conditions of 
 Definition \ref{def:aOka} are fulfilled for the morphism 
 \[F=\ts\circ\phi=s\circ f\circ\phi\colon Y\times\A^m\to X.\]
 Thus, $X$ is aOka.
 
 (2)$\Rightarrow$(1).  Suppose $X$ is aOka. 
 Let $Y\subset X$ be an affine dense open subset,
 $f\colon Y\to X$ be the identical embedding,
 and $F\colon Y\times\A^m\to X$ be the extension of $f$
 that satisfies the conditions of Definition \ref{def:aOka}.  
Then $(E,p,s)$ with $E=Y\times\A^m$, 
$p={\rm pr}_1\colon E\to Y$, and $s=F$ 
is a dominating spray over $Y$. It follows that 
$X$ is locally Gromov elliptic, 
hence  Gromov elliptic by Theorem \ref{thm:equiv}.
\end{proof}
From Theorems \ref{thm:loc-stab-flex} and \ref{thm:LT}
we deduce the following corollary,
cf. \cite[Theorem~6.2]{For23}.
\begin{cor}
Every locally stably flexible smooth algebraic variety $X$ 
over $\kk$ is aOka.
\end{cor}
The following approximation result concerns a 
(not necessary complete) 
homogeneous spaces; cf. (3) 
in Theorem \ref{thm:LT}. 
\begin{thm}[\rm{Bochnak--Kucharz  
\cite[Theorem~1.1 and Corollaries~1.2-1.3]{BoKu23}}]\label{thm:BK}  $\,$

\noindent
Let $X=G/H$ be a homogeneous space of 
a linear complex algebraic group $G$,
$Y$ be a complex affine algebraic manifold, and
$K$ be a holomorphically convex  compact set in $Y$. 
Given a holomorphic map $f \colon K \to X$, 
the following conditions are equivalent:
\begin{enumerate}\item
$f$ can be uniformly approximated 
by regular maps $K\to X$. 
\item $f$ is homotopic to a regular map $K \to X$.
\end{enumerate}
In particular,  every null homotopic holomorphic map 
$K \to X$ 
can be approximated by regular maps  $K \to X$.
\end{thm}
By a \emph{regular map}  $K \to X$ one means 
the restriction to $K$ of a morphism 
$U\to X$, where $U$ is a Zariski open 
neighborhood of $K$ in $X$.

It is worth mentioning also the following 
approximation theorem. 
\begin{thm}[\rm{Demailly--Lempert--Shiffman 
\cite[Theorem~1.1]{DLS94}}]\label{thm:DLS} $\,$

\noindent Let $X$ and $Y$ be smooth algebraic varieties, 
where $Y$ is affine and $X$ is quasiprojective.
Let $D\subset X$ be a Runge domain, 
i.e. every holomorphic
function on $D$ can be approximated, 
uniformly on compacts in $D$,  by 
holomorphic functions on $Y$. 
Then every holomorphic map
$D\to X$ can be approximated, 
uniformly on compacts in $D$, 
by Nash algebraic maps. 
\end{thm}
Recall that a holomorphic map $f\colon U\to X$ 
from an open domain $U\subset Y$
is \emph{Nash algebraic} if its graph 
$\Gamma(f)\subset U\times X$ is contained in 
a closed algebraic subset $Z\subset Y\times X$ 
of dimension $\dim(Z)=\dim(Y)$.
In fact, the approximation in 
Theorem \ref{thm:DLS} can be 
accompanied by an interpolation on 
a fixed submanifold of $Y$, 
see \cite[Theorem~1.1]{DLS94}. 

Finally, we address Kusakabe \cite{Kus21a}
for a complex algebraic version of 
Thom's jet transversality theorem and its applications. 
\subsection{Domination by affine spaces}
According to Theorems \ref{thm:Forst} 
and \ref{thm:Kusak} below, 
Gromov elliptic 
algebraic manifolds 
are dominated by affine spaces. 
We use the following notation. 
Given a surjective morphism $f\colon Y\to X$ 
of smooth algebraic varieties, 
we let $D_{\rm smooth}(f)$ stand 
for the subset of points $y\in Y$ such that 
$df|_{T_yY}\colon T_yY\to T_{f(y)}X$ 
is onto. If $\dim(X)=\dim(Y)$ 
and $f$ is dominant, then 
$D_{\rm smooth}(f)$
 is the complement of the ramification divisor of $f$. 
 In general, $D_{\rm smooth}(f)$ is 
 the maximal open subset $U$ in $Y$
 such that the restriction of $f|_U$
 is a smooth morphism onto its image.
\begin{thm}[\rm{Forstneri\v{c} \cite[Theorem 1.1]
{For17a}}]\label{thm:Forst}  
Let $X$ be a complete 
smooth complex algebraic variety 
of dimension $n$.
If $X$ is Gromov elliptic, then $X$ 
admits a  morphism $f\colon\CC^n\to X$ 
such that the restriction 
$f|_{D_{\rm smooth}(f)}\colon D_{\rm smooth}(f)\to X$ is surjective.
\end{thm}
 The proof exploits the approximation provided 
 by Theorem \ref{thm:approx}, and the latter
involves transcendental tools. 
The next result
is valid over any algebraically closed field 
$\kk$ of characteristic zero.
\begin{thm}[\rm{Kusakabe \cite[Theorem~1.2]{Kus22a}}]
\label{thm:Kusak} $\,$

\noindent
Every Gromov elliptic smooth algebraic variety
$X$ of dimension $n$ admits a morphism 
$f\colon\A^{n+1}\to X$ 
such that the restriction $f|_{D_{\rm smooth}(f)}$ 
is surjective.  
\end{thm}
See also \cite[Remark 1.9.4]{AKZ24} for a 
modified and short proof 
of Theorem \ref{thm:Kusak}
 in the case of a complete variety.
This theorem immediately leads to 
the following interpolation result. 
\begin{cor}[\rm{Kusakabe \cite[Corollary~1.5]{Kus22a}}]
\label{cor:kusak} $\,$

\noindent Let $X$ be a  Gromov elliptic smooth algebraic variety, $Y$ 
be a quasiaffine algebraic variety, and $Z \subset Y$
 be a zero-dimensional subscheme. 
Then for every morphism $f \colon Z \to X$
there exists a morphism $\tilde f \colon Y \to X$ 
such that $\tilde f|_Z = f$.
\end{cor}
Using Corollary \ref{cor:kusak} we can deduce the following 
weak version of Theorem \ref{thm:high-transitivity} 
for Gromov elliptic varieties. Let $\End(X)$ stand for 
the monoid of regular self-maps $X\to X$.
\begin{cor}[\rm{Kaliman-Zaidenberg 
\cite[Proposition~6.1]{KaZa23b}}] $\,$

\noindent Let $X$ be a smooth quasiaffine algebraic variety. If $X$ 
is Gromov elliptic, then ${\rm End}(X)$ 
acts $m$-transitively 
on $X$ for every $m\ge 1$. 
\end{cor}
See also Arzhantsev \cite{Arz23},  
\cite{Arz24},  Barth \cite{Bar23} and Kusakabe 
\cite[Corollary~1.4]{Kus22a}
for examples 
of affine varieties that admit surjective 
morphisms from affine spaces. 
However,  in some of these examples 
the surjectivity of morphisms 
restricted to their smooth loci
is not guaranteed. 
 \section{Gromov ellipticity and birational geometry}
\subsection{Gromov ellipticity versus (uni)rationality}
Recall that an algebraic variety is 
called \emph{unirational} if 
it admits a dominant rational map 
from a projective space.
An elliptic algebraic variety 
$X$ is unirational. 
Indeed, let $(E,p,s)$  be a dominating spray on $X$.
Then each fiber $E_x=p^{-1}(x)$
is an affine space which dominates $X$.
Gromov  \cite[Sec. 3.5E$''$]{Gro89} asked  
whether the opposite is true:
\begin{ques}\label{ques:unirat}
\emph{Is every smooth (uni)rational 
complete algebraic variety Gromov elliptic? }
\end{ques} 
More generally, one can ask:
\begin{ques}\label{ques:rat-conn}
\emph{Is every smooth rationally connected 
complete  algebraic variety 
Gromov elliptic? }
\end{ques} 
Since any Gromov elliptic manifold is
unirational, an affirmative answer to Question 
\ref{ques:rat-conn}
would imply that a rationally 
connected algebraic variety is 
unirational, thus resolving in the affirmative
the old open problem on coincidence of the
unirationality and the rational connectedness. 
However, the answer to the latter problem is
expected to be negative.

We say that a projective variety is
\emph{special} if it does not 
admit a dominant rational map to 
a variety of general type,
 cf. Campana \cite[Definition 2.1.2]{Cam04}.
Gromov \cite[Sec. 3.4.F]{Gro89} proposed 
``the most optimistic'' conjecture:

\medskip

\noindent {\bf Conjecture.} 
\emph{Every special smooth 
projective variety is analytically elliptic.}

\medskip

See also Campana--Winkelmann \cite{CaWi15} 
for some results and conjectures 
on the relationships between 
specialness properties of Campana and Gromov ellipticity. 
\subsection{Is Gromov ellipticity birationally invariant?}
This is a question of Gromov, 
see \cite[Remark~3.5.E$'''$]{Gro89}.
More specifically, we consider  the following question. 
\begin{ques}[\rm{cf. \cite[Remark~2(f)]{LT19}}]
\emph{Is Gromov ellipticity 
a birational property in the category 
of smooth complete algebraic varieties and 
compositions of blowups and blowdowns 
with smooth centers as morphisms? }
\end{ques} 
Indeed, a birational map between smooth 
complete varieties 
can be factored in a sequence of blowups and 
blowdowns with smooth centers, see \cite{Wlo03}. 
The behavior of Gromov ellipticity 
under blowdowns 
with smooth centers remains a mystery. However,
Gromov ellipticity 
is preserved under blowups with smooth centers 
modulo certain additional assumptions.
\begin{thm}[\rm{Kaliman--Kutzschebauch--Truong 
\cite[Theorem~0.6]{KKT18}\footnote{Cf. 
also Gromov \cite[Sec.~3.5D$''$]{Gro89} 
and L\'arusson--Truong 
\cite[Main Theorem]{LT17}.}}] $\,$

\noindent Let $X$ be an algebraic manifold 
and $Z\subset X$ be 
a smooth closed 
subvariety  of codimension $\ge 2$. 
Suppose $X$ is locally  stably flexible.
Then $X$ blown up along $Z$
is Gromov elliptic. 
\end{thm}
In general, a blowup of $X$ along a smooth center 
$Z$ does not need to preserve local  stable flexibility,
even if $X$ admits
a covering by copies of $\A^n$. 
However, this is the case if for any element 
$U_i\cong\A^n$ of the latter covering, the pair
$(U_i,Z\cap U_i)$ with a nonempty 
intersection $Z\cap U_i$ is isomorphic 
to a pair $(\A^n,\A^k)$, where 
$\A^k\subset\A^n$ is a linear subspace, 
see \cite{APS14}. 

The algebraic convex approximation property 
(the aCAP, see Remark \ref{rem:aCAP})
occurs to be stable under blowups with smooth centers. 
\begin{thm}[\rm{Kusakabe 
\cite[Corollary~4.3]{Kus20a}}] \label{thm:kusak3}
Let $X$ be a
smooth complex algebraic variety
and $A \subset X$ be a smooth closed
algebraic subvariety of codimension $\ge 2$. 
Then the blowup of $X$ along $A$ 
enjoys aCAP provided that $X$ does so. 
\end{thm}
It is not known whether the algebraic resp. 
analytic ellipticity 
of a smooth algebraic variety $X$
are preserved
under a blowup with smooth center. 
However, the blowup of $X$ 
with smooth center
is analytically elliptic provided that $X$  
is algebraically Gromov elliptic, 
see Kusakabe \cite[Corollary~1.5]{Kus20a}.
\subsection{The Poincar\'e group 
of a Gromov elliptic manifold}
Let $X$ be a complete smooth complex algebraic variety. 
 If $X$ is unirational, then $\pi_1(X)=1$ by Serre's theorem, 
 see \cite{Ser59}.
 Since every Gromov elliptic manifold is unirational, 
 $\pi_1(X)=1$ provided $X$ is Gromov elliptic.
 
 For a not necessary complete Gromov elliptic manifold $X$
 the following holds.
\begin{thm}[\rm{Kusakabe \cite[Theorem~1.3]{Kus22b}, 
\cite[Theorems~3.1 and 3.3]{Kus24}}]
\label{thm:kusak1}  $\,$

\noindent Let $X$ be a smooth complex algebraic variety. 
If $X$ is Gromov elliptic, then $\pi_1(X)$ is finite and 
 the universal covering $\tilde X$ of $X$ 
 is a Gromov elliptic algebraic variety. 
 For any finite group $\Gamma$
there exists a smooth complex quasiaffine variety $X$ 
such that 
$X$ is flexible, hence Gromov elliptic, 
and $\pi_1(X)=\Gamma$. 
\end{thm}
The following question arises; 
cf. \cite[Problem~6.4.11]{For17}.
\begin{ques}\emph{
Consider a finite morphism $X\to Y$ between 
smooth complete algebraic varieties. 
Suppose $X$ is Gromov elliptic. 
Is it true that $Y$ is  Gromov elliptic?}
\end{ques}
\subsection{Gromov  ellipticity versus uniform rationality}
\begin{defi} An algebraic variety $X$ 
of dimension $n$ is called 
\emph{uniformly rational}\footnote{There are many 
other names attributed to this same property; 
see \cite[Sec.~3]{Pop20}.}
if $X$ can be covered by open sets 
isomorphic to open sets in $\A^n$.

$X$ is called \emph{stably uniformly rational} 
 if $X\times\A^k$ is uniformly rational for some $k\ge 0$.
\end{defi}
\begin{ques}[\rm{Gromov \cite[Sec.~3.5.E$'''$]{Gro89}}] $\,$

\noindent\emph{ Is every rational smooth 
algebraic variety (stably) uniformly rational? }
\end{ques}
The affirmative answer is known for 
homogeneous spaces of 
connected affine algebraic groups, 
see V.~L.~Popov \cite[Corollary~1]{Pop20}.
As we have  already mentioned, a
complete smooth spherical variety admits an open
covering by copies of affine spaces, 
see Brion--Luna--Vust \cite[Sec.~1.5, Corollaire]{BLV86}.
Furthermore, a smooth spherical variety 
is uniformly rational, 
see  V.~L.~Popov  \cite[Theorem~4]{Pop20}. 
The total space of a locally trivial fiber bundle
over  a uniformly rational base 
with a uniformly rational general fiber is 
uniformly rational, cf.  V.~L.~Popov \cite[Theorem~2]{Pop20}. 
The next Theorems \ref{thm:Pop0} and  \ref{thm:Pop} 
provide more classes of uniformly rational
varieties.
\begin{thm}[\rm{V.~L.~Popov \cite[Theorem~1]{Pop20}}]
\label{thm:Pop0}
Let $X$ be a rational algebraic variety. If $\Aut(X)$ 
acts transitively on $X$, then $X$ is uniformly rational. 
\end{thm}
\begin{thm}[\rm{V.~L.~Popov \cite[Theorem~3]{Pop20}}]
\label{thm:Pop}
Let $G$ be a connected reductive algebraic group and $X$ 
be a smooth affine $G$-variety. Assume that every 
$G$-invariant regular function on $X$ is constant, and so 
there is a unique closed $G$-orbit in $X$.  
If this orbit is rational, then 
$X$ is uniformly rational. 
\end{thm}
The uniform rationality survives successive 
blowups with smooth centers.
\begin{thm}[\rm{Bogomolov, see
\cite[ Proposition~3.5E]{Gro89}, see also
Bodn\`ar--Hauser--Schicho--Villamayor U 
\cite[Theorem~4.4]{BHSV08} and
Bogomolov--B\"{o}hning \cite[Proposition~2.6]{BB14}}]
\label{thm:unif-rat} 

\noindent The blowup of a uniformly rational variety along 
a smooth center is uniformly rational.
\end{thm}
\begin{thm}[\rm{Arzhantsev--Kaliman--Zaidenberg 
\cite[Theorem~1.3]{AKZ24}}]\label{thm:AKZ} $\,$

\noindent 
A stably uniformly rational smooth complete 
algebraic variety $X$ is Gromov elliptic. 
\end{thm}
Due to Theorem \ref{thm:unif-rat}, 
Gromov ellipticity of a uniformly rational
smooth complete algebraic variety 
 survives successive 
blowups with smooth irreducible centers. 

From Theorems  \ref{thm:Pop0} and \ref{thm:AKZ}
we deduce such a corollary.
\begin{cor} A (stably) locally flexible smooth 
rational algebraic variety
$X$ is (stably) uniformly rational. 
If such a variety $X$ is complete, then it is Gromov elliptic. 
\end{cor}
The latter conclusion also holds for 
complete smooth $G$-varieties
verifying locally the assumptions of Theorem \ref{thm:Pop}.
\subsection{Unirationality versus  uniform rationality}
\begin{prop}[\rm{Arzhantsev--Kaliman--Zaidenberg 
\cite[Theorem~1.7]{AKZ24}}] \label{prop:unirat} $\,$

\noindent Let $X$ be a unirational complete 
variety of dimension $n$. 
Then there is  a surjective  morphism $\tilde X\to X$
from a uniformly rational complete 
variety $\tilde X$ of dimension $n$.
If $X$ is rational, then the morphism $\tilde X\to X$ 
can be chosen to be birational. 
\end{prop}
\begin{proof} By Chow's Lemma we 
may assume that $X$ is projective.
Choose a generically finite dominant rational map 
$h\colon\PP^n\dasharrow X$ which is birational if $X$ is rational.
By Hironaka's elimination 
of indeterminacy we have a commutative diagram 
\usetikzlibrary{matrix,arrows,decorations.pathmorphing}
     \begin{center}
        \begin{tikzpicture}[scale=2]
        
        \node at (1,1){$\tilde X$};
        \node at (0,0){$\PP^n$};
        \node at (2,0){$X$};
        \node at (1,0.2){$h$};
        \draw[->][] (0.8,0.8)--node[above=1pt]{$f$} (0.2,0.2); 
         \draw[->][] (1.1,0.8)--node[above=1pt]{$g$} (1.8,0.2); 
        \draw[->][thick,dashed] (0.2,0)--(1.8,0);
       
        \end{tikzpicture}
\end{center}
where $f$ is a composition of blowups with smooth 
irreducible centers. 
By Theorem \ref{thm:unif-rat}, $\tilde X$ is uniformly rational.
\end{proof}
\begin{rem}
In the case of a rational smooth projective variety $X$, 
Proposition \ref{prop:unirat} follows from 
\cite[Proposition~3.5.E"]{Gro89} due to F.~Bogomolov.
\end{rem}
\begin{cor}[\rm{\cite[Corollary~1.8]{AKZ24}}]\label{thm:unirat}
A complete algebraic variety 
$X$ is unirational 
if and only if $X$ admits a surjective morphism 
$\A^{N}\to X$ for some $N\ge\dim(X)$.
\end{cor}
\begin{proof}
By Proposition \ref{prop:unirat} 
a complete unirational $X$ is dominated 
by a complete uniformly rational variety 
$\tilde X$, where $\tilde X$ 
is Gromov elliptic due to Theorem \ref{thm:AKZ}. 
By Forstneri\v{c}' and Kusakabe's 
Theorems \ref{thm:Forst}--\ref{thm:Kusak}, 
there is a surjective 
morphism $\A^{n+1}\to\tilde X$ (resp., $\A^{n}\to\tilde X$ 
if $\tilde X$ is defined over $\CC$). 
Anyway, there is 
a surjective morphism $\A^{n+1}\to X$.
\end{proof}
We address the article \cite{BKK13} 
for a closely related subject. 
\subsection{Examples of uniformly rational algebraic manifolds} 
\begin{thm}[\rm{Bogomolov--B\"{o}hning 
\cite{BB14}}]\label{thm:BB} $\,$

(a) Every rational smooth cubic 
 hypersurface in $\PP^{n+1}$, $n\ge 2$
 is uniformly rational. 
 The same conclusion also applies to
 a small algebraic resolution of a nodal cubic 
threefold in $\PP^{4}$.

 (b)
Every smooth complete  intersection of 
two quadric hypersurfaces in 
$\PP^{n+2}$, $n\ge 3$,  is uniformly rational. 

(c) The moduli space $\overline{\mathcal{M}}_{0,n}$ 
of stable $n$-pointed rational curves
is a uniformly rational complete variety. 
\end{thm} 
Applying Theorem \ref{thm:AKZ} we deduce 
the following corollary.
\begin{cor} The varieties in Theorem \ref{thm:BB} 
are Gromov elliptic. 
\end{cor}
Cf.\ a discussion in Gromov \cite[Sec.~3.5.E$'''$]{Gro89}.
\begin{exa} A smooth cubic hypersurface 
of even dimension $n=2k$ is rational provided that
it contains a pair of skew linear $k$-spaces. 
By Theorem \ref{thm:BB}(a), 
such a hypersurface is uniformly rational,
and so Gromov elliptic by Theorem \ref{thm:AKZ}. 
See Remarks \ref{rems:cubucs} 
and the references therein 
for further examples of this type. 
Cf. also Theorem \ref{thm:cubics} below.
\end{exa}
\begin{exa}[\rm{Prokhorov--Zaidenberg \cite{PrZa23}}]
\label{ex:Fano-Mukai}
Every smooth Fano fourfold $X$ with Picard number $1$ and
of genus $10$, except at most one, up to isomorphism, 
such fourfold $X_0$, 
can be covered by copies of $\A^4$.
The exceptional fourfold $X_0$ contains 
a projective line $L$ such that
$X_0\setminus L$ is covered by copies of $\A^4$. 
Additionally, $L$ is covered in $X_0$ by open
$\A^2$-cylinders $S_i\times\A^2$,
where the $S_i$ are rational smooth affine surfaces. 
It follows that every such fourfold $X$, including $X_0$, is 
uniformly rational, and so Gromov elliptic. 
\end{exa}
\begin{cor}
Every smooth Fano fourfold $X$ 
with Picard number $1$ of
genus $10$ is Gromov elliptic. 
\end{cor}
See also Liendo--Petitjean  \cite{LP19} 
and Petitjean \cite{Pet17} for examples 
of uniformly rational affine $T$-varieties. 
\subsection{Gromov ellipticity and irrationality}
There are examples of  irrational  
smooth affine and projective
varieties that are Gromov elliptic. Let us start 
with an affine example.

\begin{exa}[\rm{cf. V.~L.~Popov 
\cite[Example 1.22]{Pop11} and \cite{Pop13}}] $\,$

\noindent Recall that an algebraic variety $X$
of dimension $n$ is called 
\emph{stably rational} if $X\times\PP^k$ is 
birationally equivalent to $\PP^{n+k}$ 
for a natural number $k$. There are irrational 
but stably rational varieties, see \cite{BCTSSD85}.
{\rm According to Saltman \cite[Theorem~3.6]{Sal84} 
(see also e.g. \cite{Bog87} and \cite{CT19})
for certain values of $n\ge 1$ and for some
finite subgroups 
$F$ of ${\rm SL}(n,\CC)$ 
the quotient $X={\rm SL}(n,\CC)/F$ 
is stably irrational (i.e. is not stably rational). 
Since $X$ is a homogeneous space of 
a semisimple algebraic group, it is flexible, 
see \cite[Proposition~5.4]{AFKKZ13}.
So, by Proposition \ref{thm:loc-flexible} $X$ is Gromov elliptic. 
Thus, $X$ is a Gromov elliptic smooth affine variety
that is stably irrational.}
\end{exa}
To deduce a projective example of this kind, it is necessary
to know the answer to the following question.
\begin{ques}\emph{
Is it true that a Gromov elliptic smooth algebraic variety 
admits a Gromov elliptic smooth completion? }
\end{ques}
The following theorem provides examples of irrational 
Gromov elliptic projective varieties, 
see Corollary \ref{cor:cubics} below.
\begin{thm}[\rm{Kaliman--Zaidenberg \cite{KaZa24b}}] \label{thm:cubics}
Every smooth cubic hypersurface $X\subset\PP^{n+1}$, $n\ge 2$, is 
Gromov elliptic. 
\end{thm}
\begin{proof}[Sketch of the proof]
Let $X\subset\PP^{n+1}$ be a smooth cubic hypersurface. 
The projection from a point $u\in X$ gives a generically 
2-to-1 rational map $\pi_u\colon X\dasharrow\PP^n$.
By permuting the pair of points on a generic fiber of $\pi_u$ 
we obtain a birational Galois involution $\tau_u$ on $X$. 

Fix $y\in X$ and choose a general $x\in X$.
The line $(xy)$ intersects $X$ 
at a third point $u$ different from $x$ and $y$. 
The involution $\tau_u$ biregularly 
sends a neighborhood $U_x$ of $x$ 
to a neighborhood $U_y$ of $y$ and sends
a line $L$ on $X$ passing through $x$ to a conic 
$C$ on $X$ passing through $y$. 

Fix $x$ and $L$ 
and vary $u'$
in a neighborhood $V$ of $u$. 
Then the image $y'=\tau_{u'}(x)$ 
runs over the neighborhood
$V'=\tau_x(V)$ of $y$. 
The image of $L$ varies in a family of conics 
$C_{y'}=\tau_{u'}(L)$ 
passing through the points
$y'=\tau_x(u')\in V'$. 

Now letting $z=L\cap \TT_uX$
and $L^*=L\setminus\{z\}$ and choosing
$x$ for zero of the vector line $L^*\cong\A^1$, 
we obtain a spray $(E,p,s)$ of rank 1
on $V'$  
which is dominating along each orbit 
$s(E_{y'})=\tau_{u'}(L^*)=C_{y'}\setminus \tau_{u'}(z)$.

We can choose $n$ lines $L_1,\ldots,L_n$ on $X$ 
passing through $x$ in $n$ independent directions.
This gives a dominating family of $n$ 
rank 1 sprays $(E_i,p_i,s_i)$ on a neighborhood 
$V_0$ of $y$ with values in $X$. 
By Gromov's Localization Lemma \ref{lem:loc}, 
each of these sprays admits an extension 
to a spray on $X$ dominating its one-dimensional orbits
inside $V_0$.

Choosing a finite covering of $X$ by 
such neighborhoods we obtain a dominating family of rank 
1 sprays on $X$. The  composition of all these 
rank 1 sprays yields a dominating spray on $X$.
Thus, $X$ is Gromov elliptic. 
\end{proof}
By Clemens--Griffiths' theorem \cite{ClGr72} 
every smooth cubic threefold 
$X$ in $\PP^4$ is irrational. 
We can therefore deduce the following corollary. 
This confirms a conjecture of Gromov \cite[Sec.~3.4.F]{Gro89}.
\begin{cor}\label{cor:cubics} 
The smooth cubic threefolds 
in $\PP^4$ are  irrational  Gromov elliptic
projective varieties. 
\end{cor}
\begin{rems}\label{rems:cubucs}
For $n\ge 3$ each smooth cubic hypersurface
 in $\PP^n$ is unirational, 
see e.g. Koll\'ar \cite{Kol02}. However, 
we do not know if a general cubic threefold
is stably irrational. 

There is a $18$-dimensional family
of rational cubic fourfolds in $\PP^5$
which contain a pair of skew planes.
Another $19$-dimensional family
of rational cubic fourfolds in $\PP^5$
consists of those cubic fourfolds which 
contain a quartic  surface scroll.
See e.g. \cite{Tre84}, \cite{Ouc20} and references 
therein for other examples. 
A general belief is that a very general  
cubic fourfold in $\PP^5$ is irrational.
However, at present, no  example 
of an irrational smooth cubic fourfold is known. 
\end{rems}
\subsection{Open questions}\label{ss:open-ques}
Likewise a smooth cubic threefold, 
a smooth quartic double solid 
admits a lot of birational involutions. 
 \begin{ques} \emph{
 Is there any Gromov elliptic smooth quartic double solid?}
\end{ques}
See \cite{CZ24} for examples of rational nodal quartic double solids
that admit Gromov elliptic algebraic small resolutions. 
 \begin{ques}[\rm{cf. \cite[Remark~2(g)]{LT19}}]
\emph{Is there a birationally 
(super)rigid Fano manifold which is Gromov elliptic? 
Especially, 
is a unirational smooth Segre quartic threefold 
in $\PP^4$ Gromov elliptic?}
\end{ques}
The superrigidity of smooth quartic threefolds in $\PP^4$
was proven by Iskovskikh and Manin, see \cite{IsMa71}; 
see also [ibid]
for Segre examples of 
unirational smooth quartic threefolds.  
 \begin{ques} \emph{
 Is the Gromov ellipticity of a smooth complete 
 algebraic variety stable under smooth deformations?
 Moreover, given a proper smooth deformation 
 family $\mathcal{X}\to S$ over a smooth base $S$,
is the locus of points $s\in S$ such that the fiber 
 $\mathcal{X}_s$ is Gromov elliptic open, 
 or closed, or constructible
 in the Zariski topology, or, over $\CC$, 
 in classical topology?
 }
 \end{ques}
 \begin{ques} \emph{
 Let $X\to S$ be a smooth morphism 
 of smooth complete varieties.
 Suppose $S$ and each fiber $X_s$, $s\in S$, 
 are  Gromov elliptic. Is it true that $X$ 
 is  Gromov elliptic?}
 \end{ques}
 The last two questions can be addressed  for (stable)
uniform rationality and (stable) local flexibility 
replacing the Gromov ellipticity.
 \section{Ellipticity of affine cones versus flexibility}
\subsection{Gromov elliptic affine cones}
\begin{thm}[\rm{Arzhantsev--Kaliman--Zaidenberg 
\cite[Theorem~1.3]{AKZ24}}]\label{thm:cones}
Let $X\subset\PP^n$ be a uniformly 
rational  projective variety and
$\hat Y\subset\A^{n+1}$ be the affine cone over $X$.
Then the punctured cone  $Y=\hat Y\setminus\{0\}$ 
is Gromov elliptic.
\end{thm}

Note that  $p\colon Y\to X$, where $p$ is 
a natural projection, 
is a principal $\GG_{\mathrm{m}}$-bundle. 
The conclusion of Theorem \ref{thm:cones} 
stays true for every
principal $\GG_{\mathrm{m}}$-bundle $Y\to X$,
provided that the associated line bundle $L\to X$
is ample or anti-ample.
However, the  assumption of ampleness is 
not necessary,
due to the next stronger result.
\begin{thm}[\rm{Kaliman \cite[Theorem~6.1]
{Kal24}}]\label{thm:kal}
Let $X$ be a complete uniformly rational algebraic variety 
and $L \to X$ be a nontrivial line bundle on $X$ 
with zero section $Z$. Then $Y = L \setminus Z$ 
is Gromov elliptic.
\end{thm}
For the trivial line bundle $L=X\times\CC$ over 
 a complex smooth complete algebraic variety $X$ we have 
$Y\cong X\times\CC^*$. So,
the group $\pi_1(Y)$ is infinite. 
Therefore, $Y$ 
cannot be Gromov elliptic, because otherwise $\pi_1(Y)$ 
must be finite by Kusakabe's Theorem \ref{thm:kusak1}. 
Thus, the hypothesis of Theorem 
\ref{thm:kal} that $L$ is nontrivial
is necessary.

There is the following analogue of 
Corollary \ref{thm:unirat} for affine cones
(which are not complete!).
\begin{thm}[\rm{Arzhantsev \cite[Theorem~1]{Arz24}}] $\,$

\noindent The affine cone $\hat Y \subset \A^{n+1}$ 
over a projective variety
$X\subset\PP^n$ admits a surjective 
morphism $\A^m\to\hat Y$ 
for some positive integer $m$ if and only if 
$\hat Y$ is unirational 
or, equivalently, $X$ is unirational. 
Furthermore, we can take $m = \dim(X) + 2$.
\end{thm}
\subsection{Flexible affine cones}
The following results certify that the flexibility 
of punctured affine cones is much stronger 
property than Gromov ellipticity. 
\begin{thm}[\rm{Kishimoto--Prokhorov--Zaidenberg 
\cite{KPZ13}, Perepechko \cite{Per13}, 
Cheltsov--Park--Won \cite{CPW16}; 
cf.~Park--Won \cite{PaWo16}}]
\label{thm:flex-cones} $\,$

\noindent Let $X=X_d$ be a del Pezzo 
surface of degree $d$  
pluri-anticanonically embedded in 
$\PP^{n}$, and let 
$Y\subset \A^{n+1}$ be the punctured affine 
cone over $X$. Then 
\begin{enumerate}
\item $Y$ is flexible for $4\le d\le 9$; 
\item $Y$ admits no nontrivial 
${\mathbb G}_{\mathrm{a}}$-action for $d\le 3$.
In particular, $Y$ is not flexible for $d\le 3$;
\item $Y$ is Gromov elliptic for every $d=1,\ldots,9$. 
\end{enumerate}
\end{thm}
The last statement follows from 
Theorem \ref{thm:cones} due to 
uniform rationality 
of the smooth complete rational surface 
$X=X_d$, see Example (7) in Section \ref{ss:flex-ex}.

The following example answers in 
negative the question in
\cite[Question~2.22]{FlZa03}.
\begin{exa}[\rm{Cheltsov--Park--Won 
\cite[Corollary~1.8]{CPW16}, 
Freudenburg--Moser-Jauslin 
\cite[Theorem~8.1(c)]{FrMJ13}}]  $\,$

\noindent Let $\bar Y\subset\A^4$ be the Fermat cubic cone
\[\bar Y=\{x_1^3+x_2^3+x_3^3+x_4^3=0\},\] 
i.e. the affine cone over the Fermat cubic surface in $\PP^3$.
Then the coordinate ring $\mathcal{O}_{\bar Y}(\bar Y)$ 
is rigid, i.e. it does not
admit any locally nilpotent derivation.
In particular, the punctured cone
$Y=\bar Y\setminus \{0\}$ is not flexible.  
However, $Y$ is Gromov elliptic due to Theorem \ref{thm:cones}.
\end{exa}
 \begin{exa} More generally,
consider a Pham--Brieskorn hypersurface 
$\bar Y\subset \A^{n+1}$, $n\ge 2$,
defined by
\[x_0^{a_0} +x_1^{a_1} +\cdots +x_n^{a_n} =0
\quad\text{where}\quad 2\le a_0\le a_1\le\cdots\le a_n.\] 
There is a conjecture \cite[Conjecture 1.22]{CPPZ21},
based on \cite{KaZa00} and \cite[P.~551 
and Example~2.21]{FlZa03}, 
claiming that the coordinate ring 
$\mathcal{O}_{\bar Y}(\bar Y)$ 
is rigid if and only if $a_1\ge 3$.
This conjecture is known to hold for $n=2$ 
\cite[Lemma~4]{KaZa00} and for $n=3$
\cite[Main theorem]{ChiDu23}, see also
\cite{ChiDa20} and \cite[Theorem~4.8.3]{Chi23}. 
\end{exa}
\begin{exa}
Let $X$ be a  Fano fourfold 
with Picard number 1 of genus 10 half-anticanonically 
embedded in $\PP^{12}$, and let $Y$ be the punctured 
affine cone over $X$. Since $X$ is uniformly rational, 
see Example \ref{ex:Fano-Mukai}, $Y$ is Gromov elliptic 
by Theorem \ref{thm:cones}.
Moreover, $Y$ is flexible, see 
Prokhorov--Zaidenberg \cite{PrZa23}.
\end{exa}
See Arzhantsev--Perepechko--S{\"u}ß 
\cite{APS14} and 
Michałek--Perepechko--S\"{u}ß \cite{MPS18} 
for further examples of flexible affine cones 
and universal torsors, i.e.
the spectra of Cox rings.
 
 \medskip
 
 \noindent {\bf Acknowledgments.} 
 It is with pleasure that we thank F.~Bogomolov, S.~Kaliman 
 and Yu.~Prokhorov for sharing their 
 expertise with the author.
 \renewcommand{\MR}[1]{}

\end{document}